\documentclass[letterpaper, 10 pt, journal, twoside]{ieeetran}
\pdfoutput=1											  

\IEEEoverridecommandlockouts                              

\def\BibTeX{{\rm B\kern-.05em{\sc i\kern-.025em b}\kern-.08em
	T\kern-.1667em\lower.7ex\hbox{E}\kern-.125emX}}
\usepackage{graphics}
\usepackage{amsmath, amssymb, mathtools, amsfonts}
\usepackage{mathabx}
\usepackage{color}
\usepackage[short]{optidef}
\usepackage{mleftright}
\usepackage[]{booktabs}
\usepackage{multirow}

\makeatletter
\let\NAT@parse\undefined
\makeatother
\usepackage[hidelinks]{hyperref}

\usepackage{amsthm}
\theoremstyle{plain}
\newtheorem{theorem}{Theorem}

\theoremstyle{definition}

\newtheorem{remark}[theorem]{Remark}

\newcommand{\bmat}{\begin{bmatrix}}
\newcommand{\emat}{\end{bmatrix}}
\newcommand{\dpartial}[2]{\frac{\partial{#1}}{\partial{#2}}}

\newcommand{\ddiff}{\mathrm{d}}

\DeclareMathOperator{\tr}{Tr}
\DeclareMathOperator{\cov}{Cov}		
\DeclareMathOperator{\EXV}{\mathbb{E}}		
\DeclareMathOperator{\vecop}{vec}		

\newcommand{\R}{\mathbb{R}}

\newcommand{\calN}{\mathcal{N}}

\newcommand{\defeq}{:=}
\newcommand{\eqdef}{=:}
\newcommand{\eye}{\mathbb{I}}  

\makeatletter
\newcommand{\subalign}[1]{%
  \vcenter{%
    \Let@ \restore@math@cr \default@tag
    \baselineskip\fontdimen10 \scriptfont\tw@
    \advance\baselineskip\fontdimen12 \scriptfont\tw@
    \lineskip\thr@@\fontdimen8 \scriptfont\thr@@
    \lineskiplimit\lineskip
    \ialign{\hfil$\m@th\scriptstyle##$&$\m@th\scriptstyle{}##$\hfil\crcr
      #1\crcr
    }%
  }%
}
\makeatother

\usepackage{tikz}
\newcommand\copyrighttext{%
  \footnotesize \textcopyright 2023 IEEE. Personal use of this material is permitted.
  Permission from IEEE must be obtained for all other uses, in any current or future
  media, including reprinting/republishing this material for advertising or promotional
  purposes, creating new collective works, for resale or redistribution to servers or
  lists, or reuse of any copyrighted component of this work in other works.
  DOI: \href{https://ieeexplore.ieee.org/document/9993720}{10.1109/LCSYS.2022.3230552}}
  
\newcommand\copyrightnotice{%
\begin{tikzpicture}[remember picture,overlay]
\node[anchor=south,yshift=5pt] at (current page.south) {\fbox{\parbox{\dimexpr\textwidth-\fboxsep-\fboxrule\relax}{\copyrighttext}}};
\end{tikzpicture}%
}

\title{
	A dual-control effect preserving formulation for nonlinear output-feedback stochastic model predictive control with constraints
	}
\author{%
	Florian Messerer$^{1}$ \and Katrin Baumg\"artner$^{1}$ \and Moritz Diehl$^{1,2}$
	\thanks{$^1$ Department of Microsystems Engineering (IMTEK), University of Freiburg, Germany
	}
  \thanks{
	{\tt\small \{firstname.lastname\}@imtek.uni-freiburg.de}
   }%
	\thanks{$^2$ Department of Mathematics,	University of Freiburg, Germany}
	\thanks{This research was supported by DFG via project 424107692 on Robust MPC and by the EU via ELO-X 953348.}
}

\begin{document}

\maketitle
\thispagestyle{empty}
\pagestyle{empty}
\copyrightnotice

\begin{abstract}
	We propose a formulation for approximate constrained nonlinear output-feedback stochastic model predictive control.
	Starting from the ideal but intractable stochastic optimal control problem (OCP), which involves the optimization over output-dependent policies, we use linearization with respect to the uncertainty to derive a tractable approximation which includes knowledge of the output model. This allows us to compute the expected value for the outer functions of the OCP exactly.
	Crucially, the dual control effect is preserved by this approximation.
	In consequence, the resulting controller is aware of how the choice of inputs affects the information available in the future which in turn influences subsequent controls.
\end{abstract}

\section{Introduction}

Model predictive control (MPC) is a powerful form of feedback control, using a model of the controlled system to compute control inputs by solving optimal control problems (OCP) in real-time \cite{Rawlings2017}.
Through a sufficiently fast feedback loop it is in general able to account for disturbances and model-plant mismatch.
However, in its standard form, i.e., nominal MPC, it uses no explicit model of uncertainty and thus may not yield adequate results for highly uncertain systems.
Stochastic and robust MPC (SMPC resp. RMPC) are the two major paradigms for addressing this issue \cite{Rawlings2017, Heirung2018, Mesbah2016}, by explicitly accounting for the uncertain model predictions.
Nonetheless, they are usually based on static models of uncertainty: they do not consider the fact that by interaction with the system information is acquired which can be used to decrease uncertainty.
This is the subject of the field of dual control \cite{Feldbaum1960}, which investigates the dual purpose of control inputs: (a) exploitation of the currently available knowledge to steer the system into a desired state, (b) exploration by choosing controls based on the respective knowledge gain, to be exploited in the future.
The latter is only relevant if the control has a dual effect for this system, i.e., loosely defined, the control does not only affect the state but also the state uncertainty.
Although dual control is usually presented for the stochastic case, it can also be formulated in the robust paradigm \cite{Lucia2014e}.

The dual control effect can be captured by explicitly including in the OCP formulation that future control inputs will depend on the information gathered up to this point \cite{Feldbaum1965, Bar-Shalom1974, Mesbah2018}.
This amounts to optimizing over output-dependent policies,
which stands in contrast to planning only open-loop control trajectories, i.e., committing to specific control inputs for future time points irrespective of the new information available by then.
In consequence, the former formulation (output-feedback SMPC) is aware of how the choice of trajectory affects the information gain and how this in turn influences future control inputs while the latter (open-loop SMPC) is not.
This allows output-feedback SMPC to resolve the explore-exploit trade-off:
It will explore exactly as much as is beneficial with respect to the objective and the constraints, as opposed to falling back to a heuristic that encourages uncertainty reduction in general and thus may unnecessarily impede performance \cite{Bertsekas2005}.

Planning over feedback laws further counteracts the quickly increasing predicted uncertainty sets from which open-loop SMPC suffers.
This is also the motivation for state-feedback SMPC, which --
under the assumption that the system state is accessible -- 
reduces predicted uncertainty by planning over state-feedback laws \cite{Goulart2006}. These can be precomputed \cite{Mayne2011} or optimized \cite{Messerer2021, Nagy2004}.
However, in situations with significant state estimation uncertainty, this assumption is clearly violated, rendering the resulting predictions questionable.

While the output-feedback SMPC problem carries all the desired properties, it is in general intractable to solve and needs to be approximated.
Here, one distinguishes between two types of approximations \cite{Filatov2000}: implicit dual control, which uses approximations that qualitatively preserve the dual control effect, and explicit dual control, in which the dual control effect is lost, but heuristic measures to ensure exploration are taken.
Recent overviews are given in \cite{Mesbah2018, Heirung2018a}.
Closely related is also perception-aware MPC, which includes some form of knowledge of the output model in the OCP: this can be ad-hoc heuristic terms in a nominal OCP \cite{Falanga2018}, or the propagation of an estimator model alongside the nominal trajectory \cite{Platt2010, Bonzanini2021, Bonzanini2022}.

In this paper we present an approximation to the output-feedback SMPC problem for constrained nonlinear systems.
This approximation is based on linearization but, crucially, retains the dual control effect.
After linearization and for specific choices of the outer components of the OCP, i.e., stage and terminal cost as well as constraint penalization, their expected value can be computed exactly.
The linearization-based uncertainty approximation can be seen as an extension of the open-loop or state-feedback formulations used in, e.g., \cite{Messerer2021, Zanelli2021, Feng2020, Hewing2020}, to the case of output feedback.
In \cite{Platt2010, Bonzanini2021} estimation uncertainty is propagated in the OCP, but without explicitly considering a state estimation based feedback law.
An extension to \cite{Bonzanini2021}, proposed in \cite{Bonzanini2022}, encodes a form of implicit feedback via a scenario tree.
Further, a similar formulation to plan over state estimate based feedback is employed in \cite{Farina2015}, but for an additively perturbed linear system such that dual control is not a topic.

We formulate the ideal output-feedback stochastic OCP in Section~\ref{sec:ofsocp}.
In Section~\ref{sec:approx} we derive a tractable approximation, the properties of which are demonstrated on an illustrative example in Section~\ref{sec:examples}, followed by a concluding Section~\ref{sec:conclusions}.

\section{Output-feedback stochastic optimal control}
\label{sec:ofsocp}
Consider the following stochastic nonlinear system,
\begin{subequations} \label{eq:stoch_nl_sys}
		\begin{align}
			x_0 &= p_0(\xi_0),\\
			x_{k+1} &= f_{k}(x_k, u_k, w_{k}), &&k=0, \dots, N-1,\\
			y_{k+1} &= g_{k+1}(x_{k+1}, v_{k+1}), && k=0, \dots, N-1,
			\end{align}
\end{subequations}
with initial state uncertainty $\xi_0$, process noise $w = (w_0, \dots, w_{N-1})$ and output noise $v =  (v_1, \dots, v_{N})$, collectively distributed as $\xi = (\xi_0, w, v) \sim \calN(0, \eye)$.
For choosing control $u_k$ at time $k$, we have the information $I_k$ available, defined recursively as
\begin{equation}
  I_1 = (u_0, y_1), \; I_k = (I_{k-1},  u_{k-1}, y_k), \; k=2,\dots,N.
\end{equation} 
Our aim is to find a policy  $\pi = (\bar u_0, \pi_1(\cdot), \dots, \pi_{N-1}(\cdot))$ defining the controls $u_0=\bar u_0$, $u_k = \pi_k(I_k)$, $k=1, \dots, N-1$, as a function of the available information.
	Information about the initial state distribution is assumed to be implicitly given via $p_0(\xi_0)$.
	Further, we assume that no more information will become available before we have to commit to a specific $u_0$, i.e., the measurement $y_0$ is already subsumed in $p_0(\xi_0)$.
	Therefore, the decision for $u_0$ is a real valued vector $\bar u_0 \in \R^{n_u}$.

Denote the simulation of system \eqref{eq:stoch_nl_sys} under policy $\pi$ and random variable realization $\xi$ by
\begin{subequations} \label{eq:dyn_nl_ofb}
	\begin{align}
		x_0^\pi(\xi) &=  p_0(\xi_0), \;
		u_0^\pi(\xi) = \bar u_0,\\
		x_{k+1}^\pi(\xi) &= f_k(x_k^\pi(\xi), u_k^\pi(\xi), w_k),
	\end{align}
	for $k=0,\dots,N-1,$ and, for $k=1,\dots,N-1$,
		\begin{align}
		y_{k}^\pi(\xi) &= g_k(x_{k}^\pi(\xi), v_{k}),\\
		I_{k}^\pi(\xi) &= (I_{k-1}^\pi(\xi), u_{k-1}^\pi(\xi), y_{k}^\pi(\xi)  ),\\
		u_{k}^\pi(\xi) &= \pi_k(I_k^\pi(\xi)).
	\end{align}
\end{subequations}
For simplicity of notation, each of the above quantities is defined as a function of the full vector $\xi$. However, from the recursion it is clear that the resulting system is causal, i.e., no quantity depends on the future.
The policy $\pi$ should minimize some deterministic measure of the stochastic cost
\begin{equation} \label{eq:cost_stoch}
	J^\pi(\xi) =  \sum_{k=0}^{N-1} l_k(x_k^\pi(\xi), u_k^\pi(\xi)) + l_N(x_N^\pi(\xi)), 
\end{equation}
with $l_0(\cdot, \cdot), \dots, l_N(\cdot)$ convex linear quadratic,
while respecting the constraints
\begin{subequations} \label{eq:constraints}
	\begin{equation}
		h_k(x_k^\pi(\xi), u_k^\pi(\xi)) \leq 0, \quad h_N(x_N^\pi(\xi)) \leq 0,
	\end{equation}
\end{subequations}
with $h_k\colon\R^{n_x} \times \R^{n_u} \to \R^{n_{h_k}}$, for $k=0,\dots,N{-}1$, and $h_N\colon\R^{n_x} \to \R^{n_{h_N}}$.
We use $h^\pi(\xi)$ to denote their concatenation in stagewise order, with $n_h$ the dimension of the corresponding output vector.
As the disturbances may be arbitrarily large, strict constraint satisfaction for every possible disturbance can in general not be enforced.
This leads to a trade-off with the performance as measured by \eqref{eq:cost_stoch}, since a higher desire to avoid constraint violation usually leads to more conservative control inputs.
Here, we choose a linear penalty on their violation,
\begin{equation} \label{eq:constraint_penalty}
	\phi_i( \eta ) \defeq \rho_i \max(0, \eta),
\end{equation}
with constraint function value $\eta\in\R$ and weight $\rho_i$.
A popular alternative would be chance-constraints, which specify that the probability of constraint violation should be below some prespecified value \cite{Mesbah2016}.

Finally, choosing the expected value as our criterion to obtain a deterministic (and thus well-defined) objective, we obtain the output-feedback stochastic OCP
\begin{mini}
	{\pi(\cdot)}
	{ \EXV_\xi \left \{ J^\pi(\xi) + \sum_{i=0}^{n_h} \phi_i( h^\pi_i(\xi)  )    \right \}  }
	{\label{eq:ofsocp}}
	{}
\end{mini}
as the problem we would ideally like to solve.
Since this problem fully retains the dual control effect of model \eqref{eq:stoch_nl_sys}, i.e., how future information can be affected as well as the possible reactions to it, it perfectly encodes the explore-exploit tradeoff as induced by the objective.
However, due to the nonlinear transformations of the random variable $\xi$ and the optimization over general policies in infinite dimensional function spaces, this problem is in general intractable.

\begin{remark}
Here, we formulate the problem explicitly as optimization over policies, as this is in line with the approximation which we will derive in the following.
A common alternative formulation would be to define the solution via a dynamic programming recursion such that the policy is implicitly represented,
where the state is either augmented by the information vector \cite{Bertsekas2005}, or propagated as a hyperstate representing the posterior distribution after having processed the new information through a Bayesian update \cite{Feldbaum1965, Mesbah2018}.
\end{remark}
\begin{remark}
	Since the first applied control $\bar u_0$ is deterministic, it is straightforward to enforce hard constraints on $\bar u_0$.
	For simplicity of notation, we assume here that the corresponding penalty weight is chosen sufficiently high such that strict constraint satisfaction is guaranteed.
	Further, control constraints are often a property of the system and thus impossible to violate.
	In this case, a prediction with non-zero probability of control constraint violation is clearly wrong.
	Conceptually, this can be handled by including the input saturation in the dynamics, or by restricting the space of policies to those that map only to the set of feasible controls, cf., e.g., \cite{Hokayem2012}.
	However, both of these approaches will not survive the linearizations central to this paper, so we instead strive for a (negligibly) low probability of constraint violation.
\end{remark}

\section{Tractable stochastic ocp approximation }
\label{sec:approx}
We will now derive an approximation to the stochastic OCP \eqref{eq:ofsocp} based on linearization. Crucially, the dual control effect will be preserved by these approximations.

\subsection{Policy and uncertainty propagation}
We obtain the nominal trajectory $\bar x = (\bar x_0, \ldots, \bar x_N)$ with nominal outputs $\bar y = (\bar y_1, \dots, \bar y_N)$ by simulation of the nonlinear system without disturbance, $\xi = 0$,
for an open-loop control trajectory $\bar u = (\bar u_0, \dots, \bar u_{N-1})$ as
\begin{subequations} \label{eq:dyn_nom}
	\begin{align}
		\bar x_0 &=  p_0(0), \\
		\bar x_{k+1} &= f_{k} (\bar x_k, \bar u_k, 0),&& k=0, \dots, N-1,\\
		\bar y_{k+1} &= g_{k+1}(\bar x_{k+1}, 0), && k=0, \dots, N-1.
\end{align}
\end{subequations}
The uncertainty propagation is based on a linearization of the system at this nominal trajectory, yielding the linear system
\begin{subequations} \label{eq:dyn_lin_at_nom}
	\begin{align}
		x_0 - \bar x_0=\; & \hat P^\mathrm{r}_0 \xi_0,\\
			x_{k+1} - \bar x_{k+1} =\;& \bar A_k(x_k\!-\!\bar x_k) + \bar B_k (u_k\!-\!\bar u_k) 
			+ \bar\Gamma_{k}w_{k}, 
		\\
			y_{k+1} - \bar y_{k+1} =\; & \bar C_{k+1}(x_{k+1}\!-\!\bar x_{k+1}) + \bar D_{k+1} v_{k+1}, 
	\end{align}
\end{subequations}
for $k=0,\dots, N-1$,
with
$\hat P^\mathrm{r}_0 = \dpartial{p_0}{\xi_0}(0)$,
$\bar A_k = \dpartial{f_k}{x_k}(\bar x_k, \bar u_k, 0)$,
$\bar B_k      = \dpartial{f_k}{u_k}(\bar x_k, \bar u_k, 0)$,
$\bar \Gamma_k = \dpartial{f_k}{w_k}(\bar x_k, \bar u_k, 0)$,
$\bar C_k      = \dpartial{g_k}{x_k}(\bar x_k, 0)$,
$\bar D_k      = \dpartial{g_k}{v_k}(\bar x_k, 0)$,
where the dependency of this linear system on the nominal trajectory $(\bar x, \bar u)$ is indicated by the overset bar.

With respect to this linearized system the information vectors $I_1, \dots, I_{N}$ can be perfectly summarized by state estimates $\hat x_k$ with covariances $\hat P_k$, obtained from a time-variant Kalman filter.
A few manipulations of the standard Kalman filter equations \cite{Anderson1979} reveal that the estimation error evolves as 
\begin{subequations} \label{eq:dyn_est_error}
	\begin{align}
		\hat x_0 - x_0 =\;& -P^\mathrm{r}_0 \xi_0, \\
			\hat x_{k+1} - x_{k+1} =\;&  (\eye - \hat K_{k+1} \bar C_{k+1}) \bar A_k (\hat x_k - x_k)
			\\
			\omit\rlap{$\displaystyle \qquad\quad
			+ ( \hat K_{k+1} \bar C_{k+1} - \eye ) \bar \Gamma_{k} w_{k} 
			+ \hat K_{k+1}\bar D_{k+1}v_{k+1},
			\nonumber $}
		\end{align}
\end{subequations}
with Kalman filter gains $\hat K_k$ and for $k=0,\dots, N-1$. The index shift on the right-hand side is due to the estimate $\hat x_{k+1}$ being dependent on the observation $y_{k+1}$.
The Kalman gains $\hat K_k$ are computed based on the linear system matrices and thus are essentially functions of the nominal trajectory $(\bar x, \bar u)$.

Restricting the space of policies to that of linear feedback $\kappa_k(\cdot)$ based on the current state estimate,
\begin{equation} \label{eq:lin_fb}
\begin{aligned}
		u_0 =&\; \bar u_0,\\
		u_k = \kappa_k(\hat x_k) \defeq&\; \bar u_k +  K_k (\hat x_k - \bar x_k)	\\
		=&\; \bar u_k +  K_k(\hat x_k - x_k) + K_k (x_k - \bar x_k),
\end{aligned}
\end{equation}
for $k=1,\dots, N-1$ and with feedback gain matrices $K_k$,
the combined dynamics of the linearized system \eqref{eq:dyn_lin_at_nom} and the estimation error \eqref{eq:dyn_est_error} under this policy $\kappa$ are given by the augmented linear system
\begin{equation} \label{eq:aug_lin_sys}
	\tilde x_0^\kappa(\xi) = 
	\bmat \hat P^{\mathrm{r}}_0  \\ -\hat P^{\mathrm{r}}_0 \emat
	\xi_0, \quad
	\tilde x_{k+1}^\kappa(\xi) = \tilde A_k^\kappa \tilde x_k^\kappa(\xi) + \tilde \Gamma_k \tilde w_k,
\end{equation}
for $k=0,\dots, N-1$, where 
\begin{subequations}
	\begin{align}
		\tilde x_k^\kappa &\defeq \bmat x_{k}^\kappa  -\bar x_k \\ \hat x_k^\kappa- x_k^\kappa  \emat, \quad
		\tilde w_k \defeq \bmat w_{k} \\ v_{k+1} \emat,\\
		\tilde A_k^\kappa &\defeq \bmat \bar A_k + \bar B_k K_k & \bar B_k K_k \\ 0 & (\eye - \hat K_{k+1} \bar C_{k+1}) \bar A_k ) \emat, \label{eq:aug_lin_sys_A} \\
		\tilde \Gamma_k &\defeq \bmat \bar \Gamma_{k} & 0 \\ (\hat K_{k+1} \bar C_{k+1} - \eye)\bar \Gamma_{k} & \hat K_{k+1} \bar  D_{k+1} \emat.
	\end{align}
\end{subequations}
The state of this system has zero mean, $\EXV_\xi \{ \tilde x_k(\xi) \} = 0$, $k=0, \dots, N$, and covariance
\begin{equation}
	\Sigma_k = \cov_\xi\{ \tilde x_k(\xi) \} = \bmat P_k & \breve P_k^\top \\ \breve P_k & \hat P_k \emat,
\end{equation}
where
$P_k = \cov_\xi\{x_k^\kappa -\bar x_k\}$ is the predicted uncertainty for deviations of the true state $x_k$ from the nominal trajectory $\bar x_k$, $\hat P_k = \cov_\xi(\hat x_k^\kappa - x_k^\kappa)$ is the predicted estimation error covariance and $\breve P_k = \EXV_\xi\{(\hat x_k^\kappa - x_k^\kappa)(x_k^\kappa-\bar x_k)^\top\}$ is their correlation.
The overall covariance $\Sigma$ propagates as
	\begin{equation} \label{eq:dyn_cov}
	\Sigma_0 = \underbrace{\bmat \hat P_0 & -\hat P_0 \\ -\hat P_0 & \hat P_0 \emat}_{\eqdef \hat \Sigma_0(\hat P_0)}, \quad
		\Sigma_{k+1} = \underbrace{\tilde A_k^\kappa \Sigma_k \tilde A_k^{\kappa\top} + \tilde\Gamma_k\tilde\Gamma_k^\top}_{\eqdef  \psi_k(\bar x_k, \bar u_k, K_k)},
	\end{equation}
for $k=0, \dots, N-1$ and with $\hat P_0 = \hat P^{\mathrm{r}}_0\hat P^{\mathrm{r}\top}_0$.

Summarizing, we approximate the evolution of the original system \eqref{eq:dyn_nl_ofb} under policy $\pi$ by (a) propagating a nonlinear nominal trajectory $(\bar x, \bar u)$ and (b) with respect to the uncertainties considering a linearization around this nominal trajectory for which we use linear feedback based on a Kalman filter estimate. Starting from the current state estimate $\hat x_0$ associated with covariance $\hat P_0$, this results in the dynamics
\begin{subequations} \label{eq:dyn_approx}
	\begin{align}
		\bar x_0 &=  \hat x_0, &
		\bar x_{k+1} &= f_{k} (\bar x_k, \bar u_k, 0),\\
		\Sigma_0 &= \hat \Sigma_0(\hat P_0), &
		\Sigma_{k+1} &= \psi_k(\bar x_k, \bar u_k, K_k),
\end{align}
\end{subequations}
for $k=0, \dots, N-1.$
The system state at time $k$ is approximately distributed as $x_k \sim \calN(\bar x_k, P_k)$, where $P_k$ is the upper left block of $\Sigma_k$.
Due to the zero block in the lower left of \eqref{eq:aug_lin_sys_A} the estimation covariance $\hat P_k$ evolves independently from the rest of the linearized system.
Via the feedback law \eqref{eq:lin_fb} the estimation covariance $\hat P_k$ is fed back into the underlying system such that larger estimation uncertainty leads to larger predictive uncertainty of the true system state $x_k$.
While for linear systems the evolution of uncertainty is independent of the state trajectory, in our case the linear system depends on the nominal trajectory $(\bar x, \bar u)$.
Therefore the approximation retains a mechanism through which (a) the choice of nominal controls $\bar u$ influences the estimation uncertainty and (b) the estimation uncertainty impacts the predictive uncertainty of the true state.
In consequence, the dual control effect of \eqref{eq:stoch_nl_sys} is qualitatively preserved. 

\subsection{Expectation of cost and constraint penalization}
We will now compute the expected value of the objective in \eqref{eq:ofsocp}, which -- with respect to the approximated dynamics \eqref{eq:dyn_approx} -- we can do exactly.
The augmented linear system \eqref{eq:aug_lin_sys} in combination with linear feedback \eqref{eq:lin_fb} follows a normal distribution for both state and control,
\begin{equation} \label{eq:dist_normal_xu}
	\underbrace{\bmat x_k \\ u_k \emat}_{\eqdef z_k} \sim \calN\Bigg( \underbrace{\bmat \bar x_k  \\ \bar u_k \emat}_{\eqdef \bar z_k}, 
	\underbrace{\bmat \eye & 0 \\ K_k & K_k \emat  \Sigma_k \bmat \eye & 0 \\ K_k & K_k \emat^\top }_{ \eqdef \tilde \Sigma_k(\Sigma_k, K_k)}
	\Bigg),
\end{equation}
for $k=0, \dots, N-1$, and $x_N \sim\calN(\bar x_N, P_N )$.
We start with the convex quadratic  stage and terminal cost functions, for which the expectation can be computed as
\begin{equation} \label{eq:stage_cost_approx}
			\EXV_{z_k\sim\calN(\bar z_k, \tilde \Sigma_k)} \{ l_k(z_k) \}
			= 
			l_k(\bar z_k) + \tfrac{1}{2} \tr( B_k \tilde \Sigma_k ),
	\end{equation}
with constant Hessian $B_k \defeq \nabla^2 l_k(\cdot)$,
for $k=0,\dots,N$, $\bar z_N\defeq \bar x_N$, $\tilde \Sigma_N \defeq P_N$,
and where in a slight overload of notation we use $l_k(z_k) = l_k(x_k, u_k)$.
Due to convexity we have $B_k\succeq 0$, such that the cost on variance $\tilde \Sigma_k$ is always nonnegative.
In consequence, the approximation of the expectation of the trajectory cost \eqref{eq:cost_stoch} is
\begin{multline} \label{eq:total_cost_approx}
	\tilde J(\bar x, \bar u, \Sigma, K) \defeq
	\sum_{k=0}^{N-1} l_k(\bar x_k, \bar u_k) 
	+ \tfrac{1}{2} \tr \left(B_k \tilde \Sigma_k (\Sigma_k, K_k)  \right)\\
	+
	l_N(\bar x_N) +\tfrac{1}{2}\tr\big(B_N \tilde\Sigma_N(\Sigma_N) \big ).
\end{multline}

\begin{figure}
	\centering
	\vspace*{5pt}
	\includegraphics[width=.8\columnwidth]{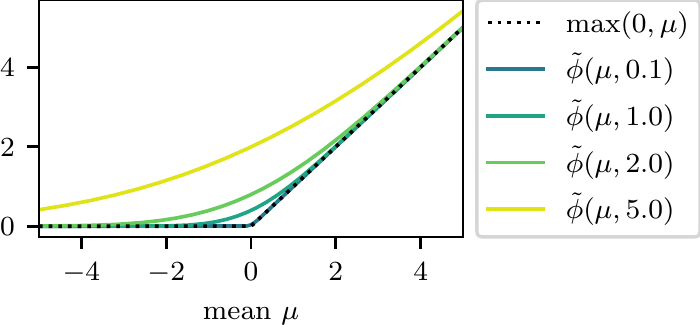}
	\caption{Illustration of the expected constraint violation $\tilde\phi(\mu, \sigma):= \EXV_{\eta\sim\calN(\mu, \sigma^2)} \{ \max(0, \eta) \}$ for several values of $\sigma$, cf. \eqref{eq:exp_max_normal}.}
	\label{fig:ev_relu}
\end{figure}

We now turn our attention to the penalization \eqref{eq:constraint_penalty} of the constraints \eqref{eq:constraints}.
An additional linearization of the constraint around the nominal trajectory leads to
\begin{equation}
\begin{aligned}
		& \EXV_{z_k\sim\calN(\bar z_k, \tilde \Sigma_k)} \{ \phi_i(h_k^i(z_k)) \} \\
		\approx &  	
		\EXV_{z_k\sim\calN(\bar z_k, \tilde \Sigma_k)} \{ \phi_i(h_k^i(\bar z_k) + \nabla h_k^i(\bar z_k)(z - z_k)  )  \} \\
		= & 	
		\EXV_{\eta \sim \calN(\bar h_k^i, \beta_k^i )} \{ \phi_i( \eta  )\}, 
	\end{aligned}
\end{equation}
for $i=1,\dots,n_{h_k}$, $k=0,\dots, N$, where $\bar h_k^i \defeq h_k^i(\bar z_k)$ is the nominal constraint value and $\beta_k^i \defeq \nabla h_k^i(\bar z_k)^\top\tilde \Sigma_k\nabla h_k^i(\bar z_k)$ the variance of the trajectory in the direction orthogonal to the constraint boundary.
For our choice of penalty function \eqref{eq:constraint_penalty} the expectation can be computed analytically as
\begin{align}
	\label{eq:exp_max_normal}
		&\EXV_{\eta\sim\calN(\mu, \sigma^2)} \{ \max(0, \eta) \} 
		= \EXV_{\nu\sim\calN(0, 1)} \{ \max(0, \sigma \nu + \mu) \}\nonumber\\
		&= \sigma \int_{-\frac{\mu}{\sigma}}^{\infty} \nu p_\calN(\nu)\ddiff \nu + \mu \int_{-\frac{\mu}{\sigma}}^{\infty} p_\calN(\nu)\ddiff \nu \nonumber\\
		&= \sigma  p_\calN\left(\frac{\mu}{\sigma}\right)+ \mu P_\calN\left(\frac{\mu}{\sigma}\right) \eqdef \tilde\phi(\mu, \sigma)
\end{align}
for $\sigma > 0$ and with $p_\calN$ the probability density function (PDF) and $P_\calN$ the cumulative distribution function (CDF) of the standard normal distribution.
Note that the expectation over the normal distribution has a smoothing effect on the originally nonsmooth penalty: the resulting function $\tilde\phi$ resembles a smoothed $\max(0, \cdot)$ in the direction of both positive $\mu$ and $\sigma$, cf. Fig.~\ref{fig:ev_relu}.
Overall this yields, as a function of the nominal value $\bar h$ and the corresponding variances $\beta$ of all constraints, the approximation to the penalty term of the stochastic OCP \eqref{eq:ofsocp} as
\begin{equation} \label{eq:constraint_penalty_approx}
	\tilde \Phi(\bar h, \beta) \defeq \sum_{k=0}^N \sum_{i=1}^{n_{h_k}} \rho_i \tilde\phi\left(\bar h_k^i, \sqrt{\beta_k^i}\right).
\end{equation} 

\begin{remark}
	Here we present only the case of quadratic stage resp. terminal cost and a linear penalty on constraint violation.
	In principle, this can be generalized to any functions for which the expectation over a normal distribution can be computed analytically, e.g., an indicator function on constraint violation (yielding a penalty on the violating probability mass)
	or a quadratic penalty.
	These functions are part of the problem definition and thus a design choice,
	and as such not discussed in this paper.
\end{remark}

\subsection{Resulting OCP formulation}
We are now ready to state the main contribution of this paper: a linearization based approximation to the original output-feedback stochastic OCP \eqref{eq:ofsocp}.
This implicit dual OCP is
	\begin{mini!}[1]
		{\bar x, \bar u, \beta,  \Sigma, K}
		{
			\begin{multlined}[t]
				\tilde J(\bar x, \bar u, \Sigma, K) + \tilde \Phi(h(\bar x, \bar u), \beta)
				\\
				+ r(K)
			\end{multlined}
		}
		{\label{eq:iocp}}
		{}%
		\addConstraint{\bar x_0}{=\hat x_0, \;\Sigma_0 = \hat\Sigma_0(\hat P_0)}
		\addConstraint{\bar x_{k+1}}{=f_k(\bar x_k, \bar u_k, 0),}{\;\;k=0,\ldots,N-1}
		\addConstraint{\Sigma_{k+1}}{ = \psi_k(\bar x_k, \bar u_k, \Sigma_k, K_k), \label{eq:iocp_unc_dyn}  }{\;\;k=0,\ldots,N-1}
		\addConstraint{ 0 }{\geq h^u_k(\bar u_k),\label{eq:iocp_constr_u_nom} }{\;\; k=0,\dots, N-1  }
		\addConstraint{\beta}{\geq  \varepsilon_{\sigma}^2 \mathbf{1} \label{eq:iocp_beta_nonneg} }
		\addConstraint{\beta_k}{\geq
		H_k(\bar x_k, \bar u_k, \Sigma_k, K_k), \label{eq:iocp_var_hk}}{\;\; k=0,\ldots,N-1 }
		\addConstraint{\beta_N}{\geq H_N(\bar x_N, \Sigma_N), \label{eq:iocp_var_hN}}
	\end{mini!}
where $\beta = (\beta_0, \dots, \beta_N)$, $\beta_k\in\R^{n_{h_k}}$, $\Sigma = (\Sigma_0, \dots, \Sigma_N)$, $K=(K_1, \dots, K_{N-1})$, $K_0=0$, $\mathbf{1} = (1,\dots, 1)$, and with $H_k$ the concatenation of the variances in constraint direction, given by, for $i=1,\dots, n_{h_k}$,  $k=0,\dots,N$,
\begin{subequations} \label{eq:constraint_var_func}
	\begin{gather}
	\begin{multlined}
			H_k^i(\bar x_k, \bar u_k, \Sigma_k, K_k) \defeq \\
			\qquad\qquad
			\nabla h_k^i(\bar x_k, \bar u_k)^\top \tilde \Sigma_k(\Sigma_k, K_k) \nabla h_k^i(\bar x_k, \bar u_k),
	\end{multlined}
	\\ 
	H_N^i(\bar x_N, \Sigma_N) \defeq 
			\nabla h_N^i(\bar x_N)^\top \tilde \Sigma_N(\Sigma_N) \nabla h_N^i(\bar x_N).
	\end{gather}
\end{subequations}
The first term in the objective is the approximation to the expectation of the trajectory cost, cf. \eqref{eq:total_cost_approx}, followed by the penalty on expected constraint violation \eqref{eq:constraint_penalty_approx}.
Here, $\beta$ is an intermediate slack variable corresponding to the variance in constraint direction.
Its purpose is to allow for a decoupling of the uncertainty dynamics 
\eqref{eq:iocp_unc_dyn} from the square-root in \eqref{eq:constraint_penalty_approx}.
While the variances $\Sigma$ are guaranteed to be positive semi-definite at every solution to \eqref{eq:iocp}, they may have arbitrarily negative eigenvalues throughout the solver iterations, since the dynamics \eqref{eq:iocp_unc_dyn} are necessarily feasible only at solutions.
In consequence, \eqref{eq:constraint_var_func} may output arbitrarily negative values which would be passed on to the square-root in \eqref{eq:constraint_penalty_approx}.
Nonnegativity of the slack variable $\beta$ \eqref{eq:iocp_beta_nonneg} on the other hand can be easily enforced throughout the iterations \cite{Nocedal2006}.
Here, the regularization parameter $\sigma_\varepsilon^2$ is included for numerical robustness, and may be seen as a minimum variance to be considered for each constraint.
Since \eqref{eq:exp_max_normal} resp. \eqref{eq:constraint_penalty_approx} impose a larger penalty for larger variance, this acts as downward pressure on $\beta$.
Further, a regularization term for discouraging aggressive feedback, 
\begin{equation} \label{eq:iocp_regularization_term}
		r(K) = \varepsilon_K \sum\nolimits_{k=1}^{N-1} \lVert K_k \rVert_\mathrm{F}^2,
\end{equation}
was found to be beneficial for convergence, but the weight should be chosen small as to not perturb the solution too strongly.
Finally, $h_k^u(\cdot)$ in \eqref{eq:iocp_constr_u_nom} denotes the constraints which are purely on the controls,  enforced on the nominal trajectory.
This can be seen as guiding constraints for the solver, since -- for sufficiently large constraint penalization weight $\rho_i$ --  these constraints are already ensured by the penalty function $\tilde \Phi$ such that \eqref{eq:iocp_constr_u_nom} will be inactive or only weakly active \cite{Nocedal2006}.

In summary,
we essentially used two types of simplification
to obtain a tractable approximation of \eqref{eq:ofsocp}:
(i) restriction of the decision space to that of affine feedback of the given structure, which -- after elimination of the Kalman gains -- is parametrized by $\bar u$ and $K$.
This results in a suboptimal solution compared to the original problem but does not introduce error per se.
(ii) linearization of the dynamics with respect to the uncertainties in order to compute the expectation (while keeping the outer structure of cost and constraint penalty intact).
This neglects the curvature of the dynamics and will be less valid for more nonlinear dynamics.
Further, the linearization is computed at the approximated mean but the full distribution is propagated with respect to this linearization.
The higher the variance, i.e., the higher the expected distance of a sample from the mean, the less valid will this linearization be. 

\begin{remark}
	Note that the variance regularization $\varepsilon_\sigma^2 > 0$ in \eqref{eq:iocp_beta_nonneg}, the feedback regularization \eqref{eq:iocp_regularization_term} as well as the guiding constraints \eqref{eq:iocp_constr_u_nom} are added for the purpose of numerical robustness of the resulting OCP. Conceptually, they are not necessary components of our proposed approximation of \eqref{eq:ofsocp}.
\end{remark}

\begin{remark}
	In the derivation of the approximation we simply posited that the Kalman filter is the adequate filter to be used with respect to the linearized system.
	While this is clear from an information perspective, it is not immediately obvious that the Kalman filter is the optimal linear filter of this structure \eqref{eq:dyn_est_error}, i.e., the structure of a Luenberger observer, to be used in the context of \eqref{eq:iocp}, as opposed to leaving the filter gain as a degree of freedom to be optimized.
	However, note that the Kalman filter achieves minimal covariance in a matrix sense \cite{Anderson1979}: If the Kalman gain $\hat K$ results in estimation error covariance $\hat P$, it holds that $\hat P \preceq \hat P'$ for all $\hat P'$ that can be achieved by any other gain $\hat K'$.
	Given that all Hessians $B_k$ in \eqref{eq:total_cost_approx} are positive semi-definite, there is no advantage in having larger estimation uncertainty in \eqref{eq:iocp}. In consequence, the Kalman gain is optimal in this context.
\end{remark}

\begin{remark}
	Here we do not exploit the specific structure of \eqref{eq:iocp} arising from the separation into nominal and tube dynamics. 
	Howevever,
	introducing the augmented state $\breve x_k = (\bar x_k, \vecop \Sigma_k) \in \R^{n_x + \frac{1}{2} n_x(n_x+1)}$, $k=0,\dots,N$, where the $\vecop$ operation exploits the symmetry of $\Sigma$, and the augmented controls $\breve u_0 = \bar u_0$, $\breve u_k = (\bar u_k, \vecop K_k, \beta_k) \in \R^{n_u + n_x n_u + n_{h_k}}$, $k=1,\dots,N-1$, $\breve u_N = \beta_N$, we see that \eqref{eq:iocp} has the structure of a standard OCP in simultaneous formulation. The complexity depends on the choice of algorithm, and standard considerations apply, cf., e.g., \cite{Rawlings2017}.
	The augmented state dimension depends quadratically on the original state, which has a strong effect on the computation times already for moderate state dimensions. For OCP formulations similar to \eqref{eq:iocp} tailored algorithms without this quadratic dependency exist \cite{ Messerer2021, Zanelli2021, Feng2020}, but for the specific structure of \eqref{eq:iocp} this is reserved for future work.
\end{remark}

\section{Illustrative Example}
\label{sec:examples}
We now illustrate the properties of the proposed formulation on a simple example.
All optimization problems are modeled via the Python interface of CasADi \cite{Andersson2019}, which allows for efficient derivative computation via algorithmic differentiation, and solved with IPOPT \cite{Waechter2006}.
The code can be found at \texttt{\url{https://github.com/fmesserer/ofsmpc}}.

We consider a nonholonomic robot with position $(r^\mathrm{x}, r^\mathrm{y}) $ that can move only in direction of its current orientation $\theta$. The controls are the speed $v$ in this direction as well as its angular velocity $\omega$, 
\begin{subequations}
	\begin{equation}
		x = \bmat r^\mathrm{x} \\ r^\mathrm{y} \\ \theta \emat, \quad
		\dot x = \bmat v \cos\theta \\ v \sin\theta \\ \omega  \emat + w, \quad
		\quad u = \bmat v \\ \omega \emat.
	\end{equation}
	The noise $w$ is piecewise constant on the discretization intervals, $w(t) = w_k$ for $t\in[t_k,t_{k+1}]$ and $w_k \sim \calN(0, \Sigma_w)$.
	The full state can be measured, but crucially the measurement noise increases approximately linearly with the distance to the $r^\mathrm{x}$-axis,
	\begin{align}
		y_k &= x_k + \sigma_\mathrm{y}(x) v_k, \\\text{with}\; \sigma_\mathrm{y}(x) &\defeq 1 + 10 \left(\sqrt{ r^\mathrm{y}_k + \varepsilon^2 } - \varepsilon\right)
	\end{align}
\end{subequations}
and $v_k \sim \calN(0, \Sigma_v)$.
The aim is to drive the state as much to the left as possible,
\begin{equation}
		l_k(x_k,u_k) = r_k^\mathrm{x} + \varepsilon_u \lVert u \rVert_2^2,
		\qquad
		l_N(x_N) = r^\mathrm{x}_N,
\end{equation}
while not crossing the $r^\mathrm{y}$-axis and with banded controls,
\begin{equation}
	0 \leq r_{k+1}^\mathrm{x}, \quad
	-u_{\mathrm{max}} \leq u_k \leq  u_{\mathrm{max}}, \; k=0,\dots, N-1,
\end{equation}
with constraint violation weighted by $\rho=10^3$.
The control regularization term with $\varepsilon_u = 10^{-6}$ is added to always ensure a well-defined problem.
The best strategy will be a trade-off between moving quickly towards the left to lower the direct cost and towards the $r^\mathrm{x}$-axis to reduce estimation uncertainty.

\begin{figure}
	\centering
	\vspace*{5pt}
	\includegraphics[width=\columnwidth]{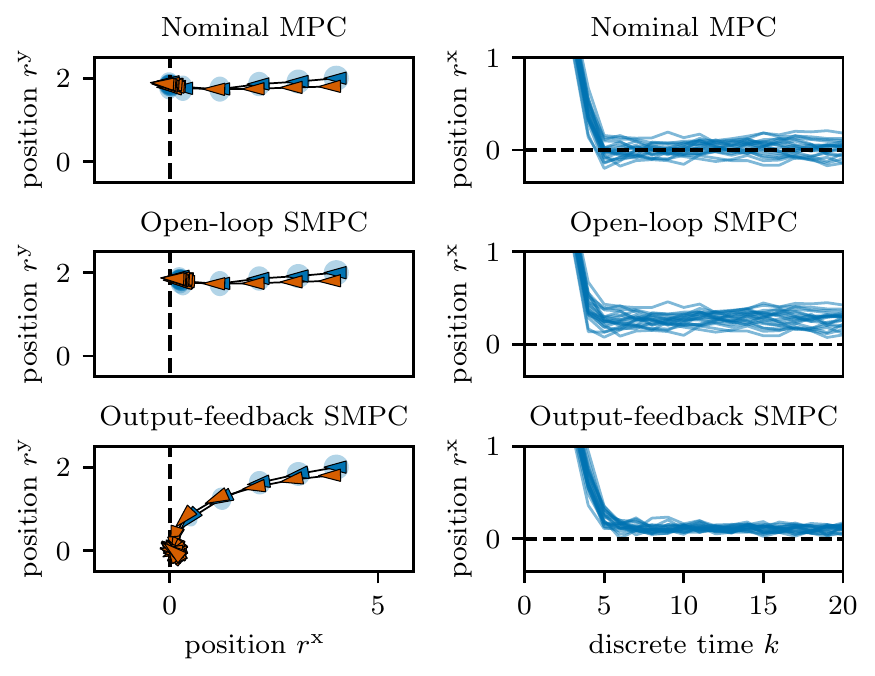}
	\caption{Left: Sample of a closed loop trajectory for each controller. The orange markers indicate the true state of the system, the blue ones the EKF estimates with the shaded regions corresponding to their 99\% confidence regions.
	Right: A closer view at the state constraint for twenty simulations per controller. Shown is the true state.
	}
	\label{fig:example_sample_trajs}
\end{figure}

We compare three MPC formulations, each with a receding horizon of $T=3$ (in continous time), discretized into $N=10$ intervals with piecewise constant controls:
(a) nominal MPC, which does not take into account uncertainty at all, (b) open-loop SMPC, which considers predictive uncertainty, but only plans an open-loop control trajectory, (c) the proposed formulation for output-feedback SMPC.
The state estimates (mean and covariance) are obtained from an Extended Kalman Filter (EKF).

Fig.~\ref{fig:example_sample_trajs} shows on its left an example trajectory of 20 timesteps for each MPC (with identical random variable realizations) and on its right takes a closer look at the state constraint over 20 sampled trajectories.
The first two controllers are not aware of the output model and thus only head left. 
Nominal MPC keeps no safety distance from the constraint, leading to frequent constraint violation.
Open-loop SMPC has no concept of actively reducing uncertainty and thus keeps a large back-off.
Only output-feedback SMPC knows the output model: it navigates to the $r^\mathrm{x}$-axis such that it is able to significantly lower the estimation uncertainty and thus the necessary backoff from the $r^\mathrm{y}$-axis.

\section{Conclusions} \label{sec:conclusions}

We presented an approximation to the ideal output-feedback stochastic OCP that preserves the dual control effect.
As the resulting formulation is aware of the output model, it can navigate to regions of the state space with low estimation uncertainty as far as is advantageous with respect to the control goal
leading to behavior qualitatively different from controllers not aware of the output model.
This is especially useful in scenarios with possibly large estimation uncertainty, in which state-feedback stochastic MPC would struggle due to a wrong underlying assumption.

\bibliographystyle{IEEEtran}
\bibliography{syscop}

\end{document}